\documentclass[a4paper,12pt,leqno]{amsart}
\usepackage{hyperref}
\usepackage{doi} 
\usepackage{url} 
\usepackage{amsmath}
\usepackage{amssymb}
\usepackage{amsthm}
\usepackage{amsopn}
\usepackage{amscd}
\usepackage{fullpage}
\usepackage[utf8]{inputenc}
\usepackage{graphics, xcolor}
\usepackage{textcomp}
\usepackage[all]{xy}
\usepackage{lscape}

\usepackage{verbatim}
\usepackage{mathrsfs}

\newcommand{\frB}{\mathfrak{B}}

\newcommand{\frM}{\mathfrak{M}}

\newcommand{\frZ}{\mathfrak{Z}}

\newcommand{\fra}{\mathfrak{a}}
\newcommand{\frb}{\mathfrak{b}}

\newcommand{\frg}{\mathfrak{g}}
\newcommand{\frh}{\mathfrak{h}}

\newcommand{\frk}{\mathfrak{k}}
\newcommand{\frl}{\mathfrak{l}}
\newcommand{\frmm}{\mathfrak{m}}
\newcommand{\frn}{\mathfrak{n}}

\newcommand{\frp}{\mathfrak{p}}

\newcommand{\frs}{\mathfrak{s}}
\newcommand{\frt}{\mathfrak{t}}

\newcommand{\bbC}{\mathbb{C}}

\newcommand{\bbG}{\mathbb{G}}

\newcommand{\bbK}{\mathbb{K}}
\newcommand{\bbL}{\mathbb{L}}

\newcommand{\bbN}{\mathbb{N}}

\newcommand{\bbR}{\mathbb{R}}

\newcommand{\bbZ}{\mathbb{Z}}

\newcommand{\caB}{\mathcal{B}}
\newcommand{\caC}{\mathcal{C}}

\newcommand{\caG}{\mathcal{G}}
\newcommand{\caH}{\mathcal{H}}
\newcommand{\caI}{\mathcal{I}}

\newcommand{\caK}{\mathcal{K}}
\newcommand{\caL}{\mathcal{L}}

\newcommand{\caP}{\mathcal{P}}

\newcommand{\caT}{\mathcal{T}}

\newcommand{\GL}{\mathbf{GL}}

\newcommand{\Std}{\mathbf{Std}}

\theoremstyle{plain}
\newtheorem{thm}{Theorem}[section]
\newtheorem{lemma}[thm]{Lemma}
\newtheorem{cor}[thm]{Corollary}
\newtheorem{prop}[thm]{Proposition}

\theoremstyle{definition}
 \newtheorem{defi}[thm]{Definition}

\theoremstyle{definition}
\newtheorem{rmks}[thm]{Remarks}
\newtheorem{rmk}[thm]{Remark}

\newtheorem{hyp}[thm]{Hypotheses}

\newcommand{\Id}{\mathrm{Id}}

\newcommand{\SL}{\mathbf{SL}}

\newcommand{\SO}{\mathbf{SO}}
\newcommand{\Sp}{\mathbf{Sp}}

\newcommand{\Ext}{\mathrm{Ext}}

\newcommand{\irr}{\mathbf{irr}}
\newcommand{\std}{\mathbf{std}}

\newcommand{\bil}[2]{\langle  #1,#2 \rangle }
\newcommand{\bilo}{\langle  .\, , . \rangle }

\def \proof {\noindent \underline{\sl Proof}. }

\begin{document}

\numberwithin{equation}{section}

\title{Generic irreducibility of parabolic induction for real reductive groups}

\author{ D. RENARD}
\address{Centre de math\'ematiques Laurent Schwartz, Ecole Polytechnique, 91128 Palaiseau Cedex, France}
\email{david.renard@polytechnique.edu}
\date{\today}

\begin{abstract} Let $G$ be a real reductive linear  group in the Harish-Chandra class. Suppose that $P$ is a parabolic
subgroup of $G$ with Langlands decomposition $P=MAN$. Let $\pi$ be an irreducible representation of the Levi factor  $L=MA$. 
We give sufficient conditions on the infinitesimal character  of $\pi$  for   the induced representation $i_P^G(\pi)$ to be irreducible.
In particular, we  prove that if $\pi_M$ is an irreducible representation of $M$, then for a
generic character $\chi_\nu$ of $A$, the induced representation $i_P^G(\pi_M\boxtimes \chi_\nu)$ is irreducible. Here the parameter $\nu$ is in 
$\fra^*=(\mathrm{Lie}(A)\otimes_\bbR \bbC)^*$ and generic means outside a countable, locally finite union of hyperplanes 
which depends only on the infinitesimal character of $\pi$.
Notice that there is no other assumption on $\pi$ or $\pi_M$ than  being irreducible, so the result is not limited
to generalised principal series or standard representations, for which the result is already well known.
\end{abstract}

\maketitle

\section{Introduction} Let $G$ be a real reductive group. We assume that there is a connected  algebraic reductive group $\bbG$ defined over $\bbR$
and that $G$ has finite index in $\bbG(\bbR)$.
Let $P$ be a parabolic subgroup of 
$G$, with Langlands decomposition $P=MAN$,  let $\pi_M$ be an irreducible representation of $M$, and $\chi$ a character of $A$.
Consider the induced representation
$\tau=i_P^G(\pi_M\otimes \chi)$ where $i_P^G$ is the functor of (normalized) parabolic induction.
\begin{thm} \label{thmintro} For generic $\chi$, the representation $\tau=i_P^G(\pi_M\otimes \chi)$ is irreducible.
\end{thm}
More precisely, since $A=\exp \fra_0$ where $\fra_0$ is a real abelian Lie algebra, characters of $A$ are of the form 
\begin{equation}\label{nu} \chi_\nu : \exp(X)\mapsto \exp(\nu(X)) ,\quad (X\in \fra_0) \end{equation}
for some $\nu\in \fra^*$. ``Generic'' in the Theorem means ``for $\chi=\chi_\nu$   with $\nu$ outside a countable, locally finite union of hyperplanes
 in $\fra^*$''. See Hypotheses \ref{hyp1} for the precise conditions when the infinitesimal character of $\pi_M$ is non singular and  Hypotheses \ref{hyp2}  
 when infinitesimal character of $\pi_M$ is singular. 
For a different perspective on potential applications, notice that  Hypotheses \ref{hyp1} and \ref{hyp2}  give in fact 
sufficient  conditions on the infinitesimal character  of  an irreducible representation $\pi$ of $L=MA$  for
   the induced representation $i_P^G(\pi)$ to be irreducible (Theorems \ref{main1} and \ref{main2}).

\medskip

 The result may seem obvious to experts, and  I was surprised  not being able to find a reference in the literature.
For $p$-adic groups,  a proof  of the analog result is given  by F. Sauvageot  in \cite{Sau}, and  a totally different one by J-F. Dat  in \cite{Dat}.
Both proofs seem difficult to adapt to the real case, however.  I propose here a very simple argument, based on a very sophisticated theory, namely 
the Kazhdan-Lusztig-Vogan theory of character multiplicities that I will  try to  describe (partly)  below. 
The motivation for writing this note came from a question of Nadir Matringe, who asked  for a reference for the result in Theorem \ref{thmintro}, 
since it is used  in his work with O. Offen and C. Yang \cite{MOY}.

After the first version of the paper was written,  David Vogan informed me that he knew 
the argument  given here.  This was  not  a surprise   since the proof consists mainly in giving  references to his work.
He also sketched   a more elementary one (in the sense that it uses less sophisticated results, on Lie algebra cohomology), but probably not 
shorter to expose from the published results. I also became aware that   the idea of using the 
  Kazhdan-Lusztig-Vogan  (KLV for short from now on) 
  algorithm to prove irreducibility  of parabolically induced representations has been used before in published works, notably 
   Matumoto \cite{MAT04} and Gan-Ichino \cite{GI}.  
I had  also thought about using this argument  in our  work with Colette Moeglin on  Arthur packets for real classical groups
\cite{MR3}. Indeed,  the last  step is the construction of arbitrary packets from packets of  ``good parity"
on a maximal Levi subgroup, by parabolic induction.  In {\sl ibid}, Thm 4.4, it is stated  that this induction  preserves irreducibility, and I 
had the vague impression that  it could be a consequence of  the ideas explained in the paper.
Eventually, we didn't use that strategy and another   quite difficult and circumvoluted   argument is given   \cite{MR6}, Thm. 5.4. 
These works prompt me to phrase  the main result  of the present paper in  Corollary \ref{coronblocks}, following the 
idea of \cite{MAT04} and  \cite{GI}.
In the final section, I go back to   \cite{MR6}, Thm. 5.4, and explain how the results here
provide some shortcuts in the  proof (and even a complete argument in most cases, 
but not more than that, some difficulties remain in some degenerate cases). 

\medskip

Let me now describe more precisely  the content of the paper.
Let $\caL$  be  a real reductive group, with Lie algebra $\frl_0$ and let $\caK$ be  maximal compact subgroup. 
Let $\frl$ be the complexified Lie algebra of $\frl_0$.
By ``representations of $\caL$",  we mean finite length Harish-Chandra modules for the pair $(\frl,\caK)$. 
Fix an infinitesimal character $\chi$ for $\caL$ and denote by $\caH\caC(\frl,\caK)_\chi$, or simply  $\caH\caC_\chi$,
 the category of representations of $\caL$ with infinitesimal character $\chi$, and by  $\bbK_\bbZ\caH\caC(\frl,\caK)_\chi$, 
 or simply $\bbK_\bbZ\caH\caC_\chi$,  its Grothendieck group with coefficients in $\bbZ$.
If $\pi$ is a representation in $\caH\caC_\chi$,  denote by $[\pi]$ its image in $\bbK_\bbZ\caH\caC_\chi$.

A result of Harish-Chandra asserts that the number of equivalence classes of irreducible representations with fixed infinitesimal character
is finite  and Langlands classification for irreducible representations, as reformulated by Vogan (see \cite{VGreen}, \cite{VIC4}), gives us a  set 
$\caP_\chi^\caL=\caP_\chi$ which parametrizes the equivalence classes  of irreducible representations in  $\caH\caC_\chi$
(in fact, it is the set ${\caP_\chi}_{/\sim\caK}$ of $\caK$-conjugacy classes in $\caP_\chi$ which is in one-to-one correspondence  
 with  equivalence classes  of irreducible representations). I will be more precise later, but for the moment, 
I will  just say that a parameter $\gamma$ in $\caP_\chi$ is roughly a character of a Cartan subgroup of $\caL$ with some additional data, from which 
one can construct  a ``standard'' representation $\std(\gamma)$ in  $\caH\caC_\chi$ (parabolically induced from a limit of discrete series modulo 
the center  of the  corresponding Levi subgroup). The standard representation $\std(\gamma)$ has a   Langlands quotient $\irr(\gamma)$
 which is irreducible, and appears with multiplicity one in $\std(\gamma)$.
Thus,    $([\irr(\gamma)])_{\gamma\in {\caP_\chi}_{/ \sim\caK} }$ is a basis of the Grothendieck group $\bbK_\bbZ\caH\caC_\chi$. 
Therefore, in the Grothendieck group, one can write for all $\delta\in \caP_\chi$,
\begin{equation} \label{matrixm} [\std(\delta)]=  \sum_{\gamma\in {\caP_\chi}_{/\sim\caK}} m(\gamma,\delta)  \, [\irr(\gamma)] \end{equation}
for some non negative integers $m(\gamma,\delta)$ (the multiplicity of $\irr(\gamma)$ in $\std(\delta)$). 

It is known from properties  of the Langlands classification relative to ``exponents''  that one can invert the linear system (\ref{matrixm}).
Since we will not use exponents in this paper,  we explain this using  the  length function $l_I$ on $\caP_\chi$ introduced by Vogan (\cite{VGreen}, 8.1.4).
Indeed,   if $m(\gamma,\delta)\neq 0$,  then $l_I(\gamma)<l_I(\delta)$, or $\gamma=\delta$,  and furthermore 
 $m(\gamma,\gamma)=1$. Therefore we can write
\begin{equation} \label{matrixM} [\irr(\delta)]=  \sum_{\gamma\in {\caP_\chi}_{/\sim\caK}} M(\gamma,\delta)  \, [\std(\gamma)] \end{equation}
for some integers $M(\gamma,\delta)$.
The Kazhdan-Lusztig-Vogan theory gives an algorithm to compute these integers $M(\gamma,\delta)$ (or equivalently the 
$m(\gamma,\delta)$). We give details about the KLV algorithm in Section \ref{KLV}.

Let us apply this   to the problem of determining when a representation $\tau=i_P^G(\pi)$  is irreducible,
for $\pi$ an irreducible representation of $L=MA$, the Levi factor of $P$, as in the beginning of this introduction.
Applying the statements in the previous  paragraph  to $\caL=L$ and to the infinitesimal character $\chi$ of $\pi$, we can write 
$\pi =\irr(\delta)$ for some parameter $\delta \in \caP_\chi^L$ and write
 \[  [\pi]= [\irr(\delta)]=  \sum_{\gamma\in {\caP_\chi^L}_{/\sim K_L}} M^L(\gamma,\delta)  \, [\std(\gamma)]. \]
By the exactness of the functor $i_P^G$, we get, 
\[  [  \tau]=[ i_P^G(\pi)]= [i_P^G(\irr(\delta))]=  \sum_{\gamma\in {\caP_\chi ^L}_{/\sim K_L}} M^L(\gamma,\delta)  \, [ i_P^G(\std(\gamma))] .\]
Now, the infinitesimal character $\chi$ for $L$ determines an infinitesimal character for $G$ that we can still denote by $\chi$.
Let us assume first that the infinitesimal character $\chi$ of $\pi$ is non singular.
Then,  a parameter $\gamma\in \caP_\chi^L$ can be extended to a parameter $\gamma^G\in \caP_\chi^G$, giving a correspondence
\begin{equation} \label{gammaLG}
 \gamma\mapsto\gamma^G, \quad 
 \caP_\chi^L \longrightarrow\caP_\chi^G,\end{equation}
so that  
\begin{equation} \label{indstdLG}  i_P^G(\std(\gamma))=\std(\gamma^G) . \end{equation}
Thus we get 
\begin{equation} \label{matrixM1}  [  \tau]=  \sum_{\gamma\in {\caP_\chi^L}_{/\sim K_L}} M^L(\gamma,\delta)  \, [ \std(\gamma^G)] \end{equation}
We can compare this to 
\begin{equation} \label{matrixM} [\irr(\delta^G)]=  \sum_{\eta\in {\caP^G_\chi}_{/\sim K}} M^G(\eta,\delta^G)  \, [\std(\eta)] \end{equation}
 to conclude that $\tau=\irr(\delta^G)$ if the following conditions are satisfied :
 \begin{itemize}
 \item[a)] The correspondence  (\ref{gammaLG})  is injective. 
 \item[b)] $M^L(\gamma,\delta)=M^G(\gamma^G,\delta^G)$ for any $\gamma,\delta\in \caP^L_\chi$.
 \item[c)] $M^G(\eta,\delta^G)=0$ if $\eta$ is not in the image of   (\ref{gammaLG}).
 \end{itemize}
We will give sufficient conditions for this to hold   (Hypotheses \ref{hyp1}).
When the infinitesimal  character $\chi$ of $\pi$ is singular, we give conditions  in  Hypotheses \ref{hyp2} so that 
the correspondence  (\ref{gammaLG}) is still well-defined and a), b), c) still hold.
 The corresponding irreducibility results are Theorems \ref{main1} and \ref{main2}, and Theorem \ref{thmintro}
is obtained as a corollary. 

\medskip

The  multiplicities $M(\gamma,\delta)$ are computed by the KLV algorithm, and this algorithm  is determined
by a set of data attached to the parameters. Under the conditions we give on the infinitesimal character,
  this set of data is preserved under the one-to-one correspondence 
$\gamma \mapsto\gamma^G$. To see this, and explain how the KLV algorithm works, we need to introduce a lot of structure theory
and results taken from Vogan's papers (about integral root systems, Cayley transforms, cross-actions, etc, in\cite{VGreen}, \cite{VIC3}, \cite{VIC4}),
 which makes the paper a little bit heavy, but the proofs consist mostly in careful bookkeeping.
 
 \medskip
 
In Section 2 and 3 of the paper, we introduce the  material to be able to describe the KLV algorithm (in case of non singular infinitesimal character).
The algorithm itself (what  the  KLV polynomials are, how they give the multiplicities  $M(\gamma,\delta)$ and how to compute them) is explained at 
the end of Section 3, and the main result for us here is that this algorithm is completely determined by the structural
data introduced in Section \ref{dataonpar}.
In section \ref{LvsG},  we show that these data are ``the same'' for $L$ and $G$, if the Hypotheses \ref{hyp1} on the infinitesimal character are satisfied.
 Indeed, we see first that the correspondence (\ref{gammaLG})
between  the Langlands-Vogan parameters  sets ${ \caP^{L^\flat}_\chi}_ {/\sim M_K^\flat} $ and  ${\caP^G_\chi}_{/\sim K}$  is  injective 
and its image is a union of blocks (this term will be explained in \S \ref{blocks}), among which is the block containing $\delta^G$. 
Then we see  that the integral root systems attached to $\chi$  are the same 
in $L$ and $G$, and likewise for all the data in Section \ref{dataonpar}.
In the last section, we show how to extend the irreducibility  result to the case of 
 singular infinitesimal character, using the ``translation data'' in Chapter 16 of \cite{ABV}. 

\medskip 
I would like to thank Nadir Matringe for initiating this work.   I would also like to thank David Vogan for useful conversations
about the problem considered here and the referees for correcting numerous either typographical  or  English grammar errors.

\section{Notation, preliminaries and structural data} \label{Not}

 For any real  Lie algebra $\frb_0$, we denote by $\frb$ its complexification. 
Let $G$ be  a real reductive group as in the introduction.
Let $\frg_0$ be the Lie algebra of $G$. We also fix a Cartan involution $\theta$ of $G$
with associated maximal compact subgroup $K$, and associated Cartan decomposition $\frg_0=\frk_0\oplus\frs_0$.
We denote by $\sigma$ the complex conjugation in $\frg$ relative to the real form $\frg_0$.
We fix a $G$-invariant non-degenerate symmetric bilinear form $\bilo$ on $\frg$ (and $\frg^*$), 
preserved by $\theta$ and which is positive definite on $\frs_0$ and definite negative on $\frk_0$.

If a group $\caG$ acts on a set $X$ and if $Y$ is a subset of $X$, we denote by $\mathrm{Centr}_\caG(Y)$ or simply $\caG^Y$ the centraliser
of $Y$ in $\caG$ and by $\mathrm{Norm}_\caG(Y)$ its normaliser and we use analogous notation for a
linear action of a Lie algebra.

\subsection{Cartan subalgebras, Cartan subgroups, roots}\label{CSG}
We recall the following well-known facts about Cartan subgroups. A Cartan subgroup $H$  is the centraliser in $G$ of a Cartan subalgebra $\frh_0$ 
of $\frg_0$.
 If the Cartan subalgebra $\frh_0$ is $\theta$-stable and decomposes as $\frh_0=\frt_0\oplus \fra_0$, then the Cartan subgroup $H$
decomposes as $H=TA$ (direct product) with $T=H\cap K$ and $A=\exp\fra_0$.
Let $\frh_0$ be a Cartan subalgebra  of $\frg_0$. Let us denote by $R(\frg,\frh)$ the root system  of $\frh$ in $\frg$, 
 by $W(\frg,\frh)$  the corresponding complex Weyl group, and by 
 $W(G,H)=\mathrm{Norm}_G(H)/H$  the real Weyl group. Depending on their values on $\frh_0$, roots are classified as real, complex or imaginary.
One can furthermore distinguish between compact imaginary  and non-compact imaginary roots (see \cite{VIC4} p. 150). Let us denote by 
\[ R_{\bbR}(\frg,\frh), \quad  R_{i\bbR}(\frg,\frh),  \quad  R_{i\bbR,c}(\frg,\frh),    \quad  R_{i\bbR,nc}(\frg,\frh),  \quad R_\bbC(\frg,\frh)\]
the subsets of $R(\frg,\frh)$ consisting of the real, imaginary,  imaginary compact,   imaginary non compact,
 and complex  roots respectively. Denote 
 by $\check \alpha=  2\frac{\alpha}{\bil{\alpha}{\alpha}}\in \frh^*$ the coroot associated to a root $\alpha\in R(\frg,\frh)$ and by $s_\alpha$
the reflection in $W(\frg,\frh)$ associated to $\alpha$.

  \medskip
  
Via Harish-Chandra isomorphisms,  any element $\lambda$ in the dual $\frh^*$ of any Cartan subalgebra $\frh$ determines
a character $\chi_\lambda$  of $\frZ(\frg)$, the center of the enveloping algebra of $\frg$ ({\sl i.e.} an infinitesimal character).
 Both $\lambda$ and $\chi_\lambda$
are said to be non singular when $\bil{\alpha}{\lambda}\neq 0$ for all $\alpha\in R(\frg,\frh)$. If a positive root system 
 $R^+(\frg,\frh)$ is given, $\lambda$ is said to be dominant if $-\bil{\check \alpha}{\lambda}\notin \bbN$ for all $\alpha\in R^+(\frg,\frh)$.

  \medskip
  
  For $\lambda\in \frh^*$, set $R(\lambda)=\{ \alpha \in R(\frg,\frh)\, \vert \,  \bil{\check \alpha}{\lambda}\in \bbZ  \}$,
the set of integral roots (for $\lambda$) and  let $W(\lambda)=W(R(\lambda))$  be the Weyl group of the root system $R(\lambda)$.
 If $\lambda$ is non singular, then put 
\[ R^+(\lambda)=\{ \alpha \in R(\lambda)\, \vert \,  \bil{\check \alpha}{\lambda}>0  \}\]
 and let $\Pi(\lambda)$ and $S(\lambda)$ be respectively the set of simple roots in $R^+(\lambda)$ and the set of simple reflections
 in $W(\lambda)$.
\medskip

In order to be able to compare roots and Weyl groups on different Cartan subalgebras, we will use the abstract Cartan subalgebra $\frh_a$ of 
$\frg$ (see \cite{VIC4}, Section 2).  We fix a positive root system $R^+(\frg,\frh_a)$ in $R(\frg,\frh_a)$ 
and a non singular dominant weight $\lambda_a \in \frh_a^*$. This defines an infinitesimal character $\chi_{\lambda_a}$. We also define
$R^a=R(\lambda_a)$, $R^{a,+}$, $W^a$, $\Pi^a$, $S^a$ to be respectively the abstract integral root system, the 
abstract integral positive root system, the abstract integral Weyl group, the abstract set of simple roots, and the abstract set of 
simple reflections.

If $\frh$ is any Cartan subalgebra of $\frg$ and $\lambda\in \frh^*$ is such that $\chi_\lambda=\chi_{\lambda_a}$, there is an isomorphism  
$i_\lambda: \frh_a^*\rightarrow \frh^*$ sending $\lambda_a$ to $\lambda$ which induces isomorphisms $R^a\rightarrow R(\lambda)$, 
$W^a\rightarrow W(\lambda)$, and so forth.

\subsection{Parabolic subgroups}\label{Psub}
Let $P$ be a parabolic subgroup of $G$ with Langlands decomposition $P=MAN$ and Levi factor $L=MA$ (direct product).
Denote by $\frp_0$, $\frmm_0$, $\fra_0$, $\frn_0$ and $\frl_0$ the respective Lie algebras of $P$, $M$, $A$, $N$ and $L$. 
We also introduce the opposite parabolic $P^-$ and its Lie algebra $\frp_0^-$.
Conjugating  with  an element of $G$, we may assume that $\frl_0$ is $\theta$-stable and $M_K:=L\cap K=M\cap K$
is a maximal compact subgroup of $L$ and $M$. Both $L$ and $M$ are in the class of groups defined in the introduction.

Let $\frh_0$ be a $\theta$-stable Cartan subalgebra of $\frl_0$. It decomposes as 
\[ \frh_0=\frh_{M,0}\oplus \fra_0= \frt_0\oplus \fra_{M,0}\oplus \fra_0.\]
 Let $R(\frn,\frh)$ be the set of roots in $R(\frg,\frh)$ 
such that the corresponding root space is in $\frn$.  Then 
\[R(\frg,\frh)=R(\frl,\frh)\coprod R(\frn,\frh) \coprod (-R(\frn,\frh)).\]
The roots $\alpha\in R(\frn,\frh)$ are either real, or complex with $\sigma(\alpha)=-\theta(\alpha)$ also in   $R(\frn,\frh)$. 
Let us choose a positive root system $R^+(\frl,\frh)$ and set $R^+(\frg,\frh)=R^+(\frl,\frh)\coprod R(\frn,\frh) $.
By \cite{KV}, Lemma 11.13 and (11.12) there is an element $h_{\delta(\frn)}\in \fra_0$ such that  
$\frl_0=\frg_0^{h_{\delta(\frn)}}$. Therefore, as $\fra_0$ is central in $\frl_0$, 
we have $\frl_0=\frg_0^{\fra_0}$.

 Since $L=MA$ it is clear that $L\subset G^{\fra_0}=G^A$. If $g\in G^{\fra_0}=G^A$, it preserves
$\frl_0$, $\frn_0$ and $\frn_0^-$ which are stable under the adjoint action of $\fra_0$. Therefore, 
since $L=\mathrm{Norm}_G(\frp) \cap \mathrm{Norm}_G(\frp^-)$, we get  
\begin{equation} \label{LAG} L=G^A=G^{\fra_0}.\end{equation}
Similar consideration apply to the complex connected group $\bbG(\bbC)$, and there we have
\begin{equation} \label{LAG2} \bbL(\bbC):=\mathrm{Norm}_{\bbG(\bbC)} (\frp) \cap \mathrm{Norm}_{\bbG(\bbC)}  (\frp^-)=\bbG(\bbC)^\fra=\bbG(\bbC)^A.
\end{equation}

Some parabolic subgroups called cuspidal are attached to Cartan subgroups : let $H=TA$ be a $\theta$-stable Cartan subgroup. Set $L=G^A$ 
 Then $L$ is a Levi factor of parabolic subgroups $P=LN$ of $G$, with Langlands decomposition $L=MA$, and $T$ is a compact Cartan subgroup
 of $M$.

\subsection{Cayley transforms}\label{SCT}
For the results in this section, we refer to \cite{VGreen}, \S 8.3.
Suppose that $H=TA$ is a $\theta$-stable Cartan subgroup of $G$, and let $\alpha\in R(\frg,\frh)$ be a real root.
Then the root vectors for $\alpha$ generate a subalgebra of $\frg_0$ isomorphic to $\frs \frl(2,\bbR)$ and we get a Lie algebra morphism
\[\phi_\alpha: \, \frs \frl(2,\bbR) \longrightarrow \frg_0\]
satisfying $\phi_\alpha(-{}^tX)=\theta(X)$.
We choose $\phi_\alpha$ so that 
\[ Z_\alpha:=\phi_\alpha  \begin{pmatrix} 1&0\\0&-1 \end{pmatrix}  \in \fra_0\subset \frh_0\quad \text{and} \quad
\phi_\alpha  \begin{pmatrix} 0&1\\0&0 \end{pmatrix}  \in \frg_0^\alpha. \]
Set
 \[ \tilde Z_\alpha:=\phi_\alpha  \begin{pmatrix} 0&1\\-1&0 \end{pmatrix}  \in \frk_0.\]
Since $G$ is linear, this map exponentiate to a group morphism :
\begin{equation} \label{Phia} \Phi_\alpha: \SL(2,\bbR) \longrightarrow G.  \end{equation}
Put 
\[ \sigma_\alpha=  \Phi_\alpha  \begin{pmatrix} 0&1\\-1&0 \end{pmatrix}, \quad m_\alpha=  \sigma_\alpha^2.\]
Then $m_\alpha\in T$ and $\sigma_\alpha\in K$. If $\alpha$ is real, then $\sigma_\alpha$ normalizes $H$ and represents  $s_\alpha\in W(G,H)$.

Define $\frh_0^\alpha=\frt_0^\alpha\oplus \fra_0^\alpha $ by setting $\fra_0^\alpha=\{X\in \fra_0\, \vert\,  \alpha(X)=0\} $
and $\frt_0^\alpha=\frt_0\oplus \bbR\widetilde Z_\alpha$. The corresponding Cartan subgroup will be denoted $c_\alpha(H)=H^\alpha=T^\alpha A^\alpha$.
Notice that $\sigma_\alpha\in T^\alpha$ and $m_\alpha\in T^\alpha\cap T$.

Let $\tilde \alpha=c_\alpha(\alpha)$ be the non compact imaginary root of $\frh^\alpha$ in $\frg$ supported on $\widetilde Z_\alpha$ (the Cayley transform of $\alpha$). 

\medskip

If  $H_1=T_1A_1$ is a $\theta$-stable Cartan subgroup of $G$, and  $\beta\in R(\frg,\frh)$ is  a non-compact imaginary  root, one can also define
a Lie algebra morphism
\[\phi_\beta: \, \frs \frl(2,\bbR) \longrightarrow \frg_0\]
which exponentiate to a group morphism 
\begin{equation} \label{Phib} \Phi_\beta: \SL(2,\bbR) \longrightarrow G,  \end{equation}
another Cartan subgroup $c^\beta(H_1)=H_1^\beta=T_1^\beta A_1^\beta$ and a real root $\tilde \beta=c^\beta(\beta)\in R(\frg,\frh_1^\beta)$.
The two constructions are inverse of each other : if  $H=TA$ is a $\theta$-stable Cartan subgroup of $G$, and $\alpha\in R(\frg,\frh)$ is  a real root
then $\Phi_\alpha=\Phi_{\tilde \alpha}$, $c^{\tilde \alpha}(c_\alpha(H))=H$ and $c^{\tilde \alpha}(\tilde \alpha)=\alpha$.
 If  $H_1=T_1A_1$ is a $\theta$-stable Cartan subgroup of $G$, and $\beta\in R(\frg,\frh_1)$ is a non compact  imaginary  root
then $\Phi_\beta=\Phi_{\tilde \beta}$, $c_{\tilde \beta}(c^\beta(H_1))=H_1$  and $c_{\tilde \beta}(\tilde \beta)=\beta$.

The following equivalent conditions define type I roots (for a real root $\alpha$ or the corresponding non compact imaginary root $\tilde \alpha$) : 
\begin{itemize} 
\item[a)] the reflection $s_{\widetilde \alpha}$ does not belong to $W(G,H^\alpha)$,
\item[b)] $T^\alpha\cap T=T$,
\item[c)] $\alpha :\, T\rightarrow \{\pm 1\} $ is not onto.
\end{itemize}

The following equivalent conditions define type II roots : 
\begin{itemize} 
\item[a)] the reflection $s_{\widetilde \alpha}$ belongs to $W(G,H^\alpha)$,
\item[b)] $T^\alpha\cap T$ has index $2$ in $T$ and  $s_{\widetilde \alpha}$ has a representative in $T\setminus T^\alpha$,
\item[c)] $\alpha :\, T\rightarrow \{\pm 1\} $ is  onto.
\end{itemize}

\medskip

\subsection{Parameters for Langlands classification} \label{Lpar} We start by fixing an infinitesimal character $\chi=\chi_{\xi_a}$
by picking an element   $\xi_a\in \fra_a^*$. We assume that 
$\xi_a $ is non singular.  Abstract integral roots, etc, defined in \S \ref{CSG} with respect to an element $\lambda_a\in 
\frh_a^*$ are now  defined with respect to this element $\xi_a$.

We recall the set of parameters  ${\caP^G_\chi}_{/\sim K}$ 
for the Langlands classification 
of  irreducible representations of $G$  with infinitesimal character $\xi_a$  (see \cite{VGreen} , \cite{VIC4}). 

\begin{defi} \label{defpar} A parameter $\gamma$ is a multiplet 
\[\gamma=(H=TA,\Gamma,\bar \gamma)\]
where $H=TA$ is a $\theta$-stable Cartan subgroup of $G$,  $\Gamma$ is a character of $H$, and $\bar \gamma\in \frh^*$, 
satisfying the following conditions $a)$, $b)$, $c)$:
\begin{itemize}
\item[a)] $\bar \gamma_{\vert \frt}\in i\frt_0^*$,  and $\bil{\alpha}{\bar \gamma}\neq 0$,  $(\forall \alpha\in R_{i\bbR}(\frg,\frh))$.
\end{itemize}
Set 
\begin{align*}
&R^+_{i\bbR}=R^+_{i\bbR}(\frg,\frh)=\{  \alpha\in R_{i\bbR}(\frg,\frh)\vert \, \bil{\alpha}{\bar \gamma}> 0 \},
&R^+_{i\bbR,c}=R^+_{i\bbR}\cap R_{i\bbR,c}(\frg,\frh) ,\\
&\rho(R^+_{i\bbR})=\frac{1}{2}\sum_{\alpha\in R^+_{i\bbR}}\alpha,  &\rho(R^+_{i\bbR,c})=\frac{1}{2}\sum_{\alpha\in R^+_{i\bbR,c}}\alpha.
\end{align*}

\begin{itemize}
\item[b)] $d\Gamma=\bar \gamma+\rho(R^+_{i\bbR})-2\rho(R^+_{i\bbR,c})$.
\item[c)] The infinitesimal character  $\chi_{\bar \gamma}$  equals $\chi$.
\end{itemize}
\end{defi}

\medskip

Attached to a parameter $\gamma=(H=TA,\Gamma,\bar \gamma)$ as above, 
 there is a standard representation $\std(\gamma)$ (see \cite{VGreen} or \cite{VIC4}); 
it may be defined by parabolic induction from a discrete series representation on a
  cuspidal parabolic subgroup  $P=MAN$ attached to the Cartan subgroup $H=TA$. 
The group $N$ is chosen so that the Langlands subquotients appear as
quotients of $\std(\gamma)$. For non singular  infinitesimal character,
this quotient is irreducible and is denoted by  $\irr(\gamma)$.

Let us denote by $\caP^G_\chi$ the set of parameter $\gamma$ as above, and by  ${\caP^G_\chi}_{/\sim K}$ the  set of 
$K$-conjugacy classes in  $\caP^G_\chi$.
The Langlands classification for non singular infinitesimal character $\chi$ is the following theorem 
(see \cite{VIC4}, Thm. 2.13 and the references given there).
\begin{thm} \label{LVclas} Suppose that $\pi$ is an irreducible representation of $G$ with non singular infinitesimal character $\chi$.
Then there is a parameter $\gamma\in \caP^G_\chi$ such that $\pi=\irr(\gamma)$.
If two parameters $\gamma_1$ and $\gamma_2$ satisfy $\pi=\irr(\gamma_1)=\irr(\gamma_2)$, then 
$\gamma_1$ and $\gamma_2$ are $K$-conjugate. 
\end{thm}

We will constantly  abuse notation by denoting the $K$-conjugacy class of an element $\gamma\in \caP^G_\chi$ also by $\gamma$, and conversely, 
for a conjugacy class  $\gamma\in {\caP^G_\chi}_{/\sim K}$, we denote again by $\gamma$ the choice of a representative in   $\caP^G_\chi$.
 Usually, this should not lead to any confusion.

 \subsection{Cayley transforms and cross-action on parameters}
 
For any $\gamma=(H=TA,\Gamma,\bar \gamma)\in {\caP^G_\chi}$  and for any $w\in W(\bar \gamma)$, 
a new element $w\times \gamma$ in ${\caP^G_\chi}$, with first component $H=TA$ is defined in \cite{VGreen}, \S 8.3.
When $\alpha\in R^+(\bar \gamma)$ is a simple root, the other components of $s_\alpha\times \gamma$ are given explicitly 
in {\sl ibib.} Lemma 8.3.2. One can use the isomorphisms $i_{\bar \gamma}$ of \S \ref{CSG} to transport this to an action
of $W^a$ on $\caP^G_\chi$ and ${\caP^G_\chi}_{/\sim K}$ (see \cite{VIC4}, Section 2).

 In Section \ref{SCT}, the Cayley transform of a $\theta$-stable Cartan subgroup $H=TA$ with respect to a real root $\alpha$
 has been recalled. In  {\sl ibib.} \S 8.3 this definition is extended to Langlands parameters.

We recall first the parity conditions on real integral roots.
 If $\gamma=(H=TA,\Gamma,\bar \gamma)\in {\caP^G_\chi}$ and if $\alpha\in R(\bar \gamma)$
 is a real root, we say that $\alpha$ satisfies the parity condition with respect to $\gamma$ if 
\begin{equation} \label{paritycondition} \Gamma(m_\alpha)=\epsilon^G_\alpha (-1)^{\bil {\check \alpha}{\bar \gamma}}. \end{equation}
Here $\epsilon^G_\alpha\in \{\pm\}$ is defined
in \cite{VGreen}, Def. 8.3.11.
 
 If $\alpha$ is a real integral root 
  satisfying the parity condition, then the Cayley transform
 $c_\alpha(\gamma)$ is defined as a subset of $\caP^G_\chi$. It is a singleton if $\alpha$ is type II,
 and we set $c_\alpha(\gamma)=\{\gamma_\alpha\}$, and  if $\alpha$ is type I, then 
$c_\alpha(\gamma)=\{\gamma_\alpha^+,\gamma_\alpha^-\}$, with $s_{\tilde \alpha}\times \gamma_\alpha^\pm= \gamma_\alpha^\mp$.
 The first component of $\gamma_\alpha$ or $\gamma_\alpha^\pm$ is $H^\alpha=c_\alpha(H)$.
 
 In the other direction, one can define also Cayley transform of a parameter 
 $\gamma=(H=TA,\Gamma,\bar \gamma)\in {\caP^G_\chi}$ 
 with respect to a non compact imaginary integral root $\beta$. 
 The Cayley transform $c^\beta(\gamma) $ is a singleton if $\alpha$ is type I, and we set
  $c^\beta(\gamma)=\{\gamma^\beta\}$, and  if $\beta$ is type II, then 
$c^\beta(\gamma)=\{\gamma^\beta_+,\gamma^\beta_-\}$, with $s_{\tilde \beta}\times \gamma^\beta_\pm= \gamma^\beta_\mp$.
 The first component of $\gamma^\beta$ or $\gamma^\beta_\pm$ is $H^\beta=c^\beta(H)$.

 \subsection{Data associated to a Langlands parameter}\label{dataonpar}
We associate to any $\gamma=(H=TA,\Gamma,\bar \gamma)\in {\caP^G_\chi}$ the following set of data :

\begin{itemize}
\item[(1)] $R(\bar \gamma)$, the integral root system defined by $\bar \gamma$ and $R^+(\bar \gamma)$, the set of  integral  positive roots. 

\item[(2)]  $R^+_{i\bbR}(\bar \gamma)= R^+(\bar \gamma)\cap R^+_{i\bbR}(\frg,\frh)$,  the set of    integral imaginary positive  roots.

\item[(3)]  $R^+_{i\bbR,c}(\bar \gamma)= R^+(\bar \gamma)\cap R_{i\bbR,c}(\frg,\frh)$,  the set of    integral imaginary compact positive roots
 and  $R^+_{i\bbR,nc}(\bar \gamma)= R^+(\bar \gamma)\cap R_{i\bbR,nc}(\frg,\frh)$, the set of    integral imaginary non  compact  positive roots.

\item[(4)]  $R^{+,I} _{i\bbR,nc}(\bar \gamma)$ the set of integral imaginary non  compact positive roots of type I
and  $R^{+,II} _{i\bbR,nc}(\bar \gamma)$ the set of integral imaginary non  compact positive roots of type II.

\item[(5)]  $R_{\bbR}(\bar \gamma)= R(\bar \gamma)\cap R_{\bbR}(\frg,\frh)$,  the set of   integral real  roots and 
$R^+_{\bbR}(\bar \gamma)= R^+(\bar \gamma)\cap R_{\bbR}(\frg,\frh)$, the set of   integral real positive  roots.

\item[(6)]  $R_{\bbR,0}(\bar \gamma)$,  the set of   integral real  roots not satisfying the parity condition for $\gamma$ and 
$R_{\bbR,1}(\bar \gamma)$ is the set of    integral real  roots  satisfying the parity condition for  $\gamma$. 

\item[(7)]  $R^I_{\bbR,1}(\bar \gamma)$,  the set of   integral real  roots  satisfying the parity condition for $\gamma$ of type I
 and  $R^{II}_{\bbR,1}(\bar \gamma)$ the set of   integral real  roots satisfying the parity condition for $\gamma$ of type II.

\item[(8)]  $R^+_{\bbC,1}(\bar \gamma)$,  the set of   integral complex positive  roots such that $\theta(\alpha)\in R^+_{\bbC}(\bar \gamma)$
and  $R^+_{\bbC,0}(\bar \gamma)$, the set of    integral positive   roots  such that  $\theta(\alpha)\notin R^+_{\bbC}(\bar \gamma)$.

\end{itemize}

\begin{rmks}\label{Rmkmin}
\begin{itemize}
 \item[a)]The  integral root system  $R(\bar \gamma)$ is $\theta$-stable (\cite{VGreen}, Lemma 8.2.5).

\item[b)]   A parameter $\gamma\in \caP^G_\chi$ is said to be minimal if $R_{\bbR,1}(\bar \gamma)$ and 
$R^+_{\bbC,0}(\bar \gamma)$ are empty. It is equivalent to the fact that for any non real simple root in $R^+(\bar \gamma)$,
$\theta(\alpha)\in R^+(\bar\gamma)$ and  no  real simple root in $R^+(\bar \gamma)$ satisfies the parity condition 
  (see \cite{VGreen}, Def. 8.6.5).
  The standard representation $\std(\gamma)$ is irreducible  if and only if $\gamma$ is minimal (\cite{VGreen}, Thm. 8.6.6).

\item[c)]As we will see in the next section, the data above for all $\gamma\in {\caP^G_\chi}_{/\sim K}$ (a finite set)
is sufficient to determine the multiplicities $M(\gamma,\delta)$ or $m(\gamma,\delta)$ for any  $\gamma,\delta \in {\caP^G_\chi}_{/\sim K}$
via the KLV algorithm.

\end{itemize}
\end{rmks}

\section{The  KLV-algorithm}\label{KLV}
We follow here \cite{VIC4}, Section 12.

\subsection{Blocks }\label{blocks}
 To describe (part of) the KLV-algorithm, we first introduce the notion of block in  ${\caP_\chi^G}_{/\sim K}$.
Blocks are equivalence classes on ${\caP^G_\chi}_{/\sim K}$ for the equivalence relation generated by 
\[ \gamma_1\sim \gamma_2 \text{ if } m(\gamma_1,\gamma_2)\neq 0,  \]
{\sl i.e.} $\irr(\gamma_1)$ occurs as a subquotient in $\std(\gamma_2)$.
This is also the equivalence relation generated by the weaker condition that  $\gamma_1\sim \gamma_2 $ if $\irr(\gamma_1)$
 and $\irr(\gamma_2)$ both occur as  some subquotient of the same standard representation $\std(\delta)$.
Another characterisation is given in \cite{VGreen}, Thm 9.2.11. Block equivalence is generated by the following relations :
 if $\gamma\in \caP^G_\chi$ and $\beta$ is 
a simple non-compact imaginary root in $R^+(\bar \gamma)$, then $\gamma\sim\gamma'$ for any $\gamma'\in c^\beta(\gamma)$, and if  
$\alpha$ is  a simple complex root in $R^+(\bar \gamma)$, then $\gamma\sim s_\alpha \times \gamma$.

Via the Langlands-Vogan  parametrisation (Thm. \ref{LVclas}),
 block equivalence gives an equivalence relation on equivalence classes of irreducible representations
 which can be characterized in terms of $\Ext$ groups, namely, it is the equivalence relation generated by $\pi_1\sim \pi_2$ if 
$\Ext^1(\pi_1,\pi_2)\neq \{0\}$ (see \cite{VGreen}, Def. 9.2.1 and Prop. 9.2.10).
A result of Casselman (see \cite{VGreen}, Cor. 9.2.24) states that if  two irreducible representations $\pi_1$ and $\pi_2$ of $G$ 
are in different blocks, then 
$\Ext^*(\pi_1,\pi_2)=0$.
 Therefore the set of parameters $\caP_\chi^G$ admits a partition in blocks 
\begin{equation}\label{blocks}{\caP^G_\chi}_{/\sim K} =\coprod_i \frB_i . \end{equation}
If two parameters $\gamma,\delta  \in \caP_\chi^G$ are in different blocks, then $M(\gamma,\delta)=0$.

Let us fix a block $\frB$ in the partition above. 
\begin{defi}\label{li}
The integral length of a parameter $\gamma=(H=TA,\Gamma,\bar \gamma)\in \caP^G_\chi$ is  
\[ l_I(\gamma)=\frac{1}{2} \left\vert \{\alpha\in R^+(\bar \gamma)\vert \, \theta(\alpha) \notin R^+(\bar \gamma) \} \right\vert +\frac{1}{2}
\dim A -c^G_0\] 
\end{defi}
\begin{rmk}
The constant $c^G_0$ may be  chosen so that $ l_I(\gamma)\in \bbN$, for all $\gamma\in \frB$, 
 but the choice of $c_0^G$ is irrelevant for the KLV algorithm
since it is always the difference  $ l_I(\gamma_1)- l_I(\gamma_2)$ between the integral length
 of two parameters  $\gamma_1,\,  \gamma_2\in \caP^G_\chi$ 
which matters.
\end{rmk}

The Hecke algebra  $\caH=\caH(W^a)=\caH(W(\xi_a))$ of the abstract Weyl group $W^a$ is defined in \cite{VIC4}, Def. 12.4.
  This is an algebra over $\bbZ[u^{\frac{1}{2}}, u^{-\frac{1}{2}}]$ generated by elements $T_w$, $w\in W^a$ with the relations
  given in {\sl ibid.}

The Hecke module of $\frB$ is the free module over $\bbZ[u^{\frac{1}{2}}, u^{-\frac{1}{2}} ]$ with basis $\{\gamma\in \frB\}$.
 Let us denote this module by $\frM(\frB)$.
The action of $\caH$ on $\frM(\frB)$ is given also in {\sl ibid}. More precisely, what is given are formulas for 
$T_s\gamma$ when $\gamma\in \frB$ and $s\in S^a$ is a  simple  reflection.
This simple reflection corresponds to a simple root $\alpha\in R^+(\bar \gamma)$ via the isomorphisms $i_{\bar \gamma}$
and the formula depends on $\alpha$.
For instance, if $\alpha$ is type II real satisfying the parity condition,  then 
\[  T_s\gamma=(u-1)\gamma-s_\alpha\times \gamma +(u-1)c_\alpha(\gamma). \]

Let  $s$ be a simple reflection in $S^a$. For any  $\gamma_1$, $\gamma_2 \in \frB$, write 
$\gamma_1\stackrel{s}{\rightarrow}\gamma_2$ if $\alpha_1$,  the corresponding simple root  in $ R^+(\bar \gamma_1)$   is

$\bullet$ either complex 
with $ \theta \alpha_1\notin  R^+(\bar \gamma_1)$ and $\gamma_2=s\times \gamma_1$ 

$\bullet$ or   real, satisfying the parity condition 
with respect to $\gamma_1$ and $\gamma_2\in c_{\alpha_1}(\gamma_1)$. 

Equivalently, if $\alpha_2$ is the corresponding simple root  
in $ R^+(\bar \gamma_2)$, then $\alpha_2$ is complex and $ \theta \alpha_2\notin  R^+(\bar \gamma_2)$ or 
$\alpha_2$ is non compact imaginary with respect to $\gamma_2$ and $\gamma_1\in c^{\alpha_2}(\gamma_2)$.

If $\gamma_1\stackrel{s}{\rightarrow}\gamma_2$, then $l_I(\gamma_2)=l_I(\gamma_1)-1$, and we have also 
\[  s\times \gamma_1\stackrel{s}{\rightarrow}s\times \gamma_2, \quad \gamma_1\stackrel{s}{\rightarrow}s \times \gamma_2,
\quad  s \times \gamma_1\stackrel{s}{\rightarrow}\gamma_2\]

In \cite{VIC4}, Def. 12.12 an order relation is defined on $\caB$ and denoted $\gamma_1\leq _r\gamma_2$. Let us  recall 
some properties of this partial order relation. 
\begin{itemize}
\item[a)] If  $\gamma_1\leq _r\gamma_2$, then $l_I(\gamma_1)\leq l_I(\gamma_2)$ and if 
$\gamma_1\leq _r\gamma_2$, and  $l_I(\gamma_1)= l_I(\gamma_2)$ then $\gamma_1=\gamma_2$.

\item[b)] If $m(\gamma_1,\gamma_2)\neq 0$, or    $M(\gamma_1,\gamma_2)\neq 0$,  then $\gamma_1\leq _r \gamma_2$. 
\end{itemize}

The next lemma is \cite{VIC4}, Lemma 12.18. It is used to set up the  induction step for computing KLV polynomials.
\begin{lemma} Suppose that $\gamma, \delta\in \frB$ and $m(\gamma, \delta)\neq 0$.
Then we can find $\delta'\in \caB$ and a simple reflection $s\in W^a$ such that  $\delta\stackrel{s}{\rightarrow}\delta'$, and for any such $s$, 
one of the following conditions is satisfied
\begin{itemize}
\item[(i)] $m(\gamma,\delta')\neq 0$.
\item[(ii)] There exists $\gamma'\in \frB$ such that  $\gamma\stackrel{s}{\rightarrow}\gamma'$ and 
 $m(\gamma',\delta')\neq 0$.
\item[(iii)] Let  $\alpha$ be  the  simple root  in $ R^+(\bar \delta)$ corresponding to $s$. Then
 $\alpha$ is real, satisfying the parity condition 
with respect to $\delta$ and $(i)$ or $(ii)$ holds with $s\times \delta'$ replacing $\delta'$.
\end{itemize}
\end{lemma}

The next ingredient in the KLV algorithm is the duality map $D$ on $\frM(\frB)$ (\cite{VIC4}, Lemma 12.14).
\begin{lemma} There is a unique $\bbZ$ linear map $D: \frM(\frB)\rightarrow \frM(\frB)$ with the following properties. 
For any $\gamma\in \frB$, write 
\[D\gamma=  u^{-l_I(\gamma)} \sum_{\phi\in \frB}  (-1)^{l_I(\gamma)-l_I(\phi)} R_{\phi\gamma} \, \phi \]
for some polynomials $R_{\phi\gamma}\in \bbZ[u^{\frac{1}{2}}, u^{-\frac{1}{2}} ]$. 
Then 
\begin{itemize}
\item[a)] $D(um)=u^{-1}D(m), \; (\forall m\in \frM(\frB))$.
\item[b)] $D((T_s+1)m)=u^{-1}(T_s+1)D(m), \; (\forall m\in \frM(\frB), \, \forall s\in S^a)$.
\item[c)] $R_{\gamma\gamma}=1, \, (\forall \gamma \in \frB)$. 
\item[d)] $R_{\phi\gamma}\neq 0  \implies \phi\leq_r\gamma, \, (\forall \gamma \in \frB)$. 
\end{itemize}
The map $D$ has the following extra properties
\begin{itemize}
\item[e)] $R_{\phi\gamma} $ is a polynomial in $u$ of degree $\leq l_I(\gamma)-l_I(\phi)$.
 \item[f)] $D^2=\Id_{\frM(\frB)}$.
 \item[g)] The specialisation  of $D$ at $u=1$ is the identity.
\end{itemize}
\end{lemma}
\bigskip 

We can finally define the KLV polynomials. This is \cite{VIC4}, Lemma 12.15.
\begin{lemma} For any $\gamma\in \frB$, there is a unique element $C_\gamma= \sum_{\phi\in \frB} P_{\phi\gamma} \phi
\in \frM(\frB)$ (with  coefficients   $P_{\phi\gamma}$ in $\bbZ[u^{\frac{1}{2}}, u^{-\frac{1}{2}} ]$) such that 
$D(C_\gamma)= u^{-l_I(\gamma)} C_\gamma$,  $P_{\gamma\gamma}=1$,  $P_{\phi\gamma}\neq 0 $ only if $ \phi\leq_r\gamma$
and if $\phi\neq \gamma$, then  $P_{\phi\gamma}$ is a polynomial in $u$ of degree $\leq \frac{1}{2}\left(  l_I(\gamma)-l_I(\phi)-1\right)$.
\end{lemma}

\bigskip 

The next result proves  the Kazhdan-Lusztig-Vogan conjecture on multiplicities. 
\begin{thm}
The integers $M(\gamma,\delta)$, $\gamma,\delta\in \frB$ are  given by 
\[ M(\gamma,\delta) =  (-1)^{l_I(\delta)-l_I(\gamma)} P_{\gamma\delta}(1). \]
\end{thm}
It is proved by Vogan, Lusztig-Vogan if the infinitesimal character is integral, and an argument  by Bernstein settle the non integral case.
See \cite{VIC4} and \cite{ABV} for a discussion of this fundamental result and references to the original papers.

\medskip

Finally, there is an algorithm which computes the KLV polynomials $P_{\gamma,\delta}$. It is described in Prop. 6.14 of \cite{VIC3}.
It starts with the fact that $P_{\delta,\delta}=1$
 for any $\delta\in  {\caP_\xi^G}_{/\sim K}$
 and that  $P_{\gamma\delta}=0$ if $\gamma\nleq_r\delta$
   in ${\caP_\xi^G}_{/\sim K}$. If $P_{\gamma',\delta'}$ is known whenever $l_I(\delta')<l_I(\delta)$ or 
$l_I(\delta')=l_I(\delta)$ and $l_I(\gamma')>l_I(\gamma)$, then there are formulae for computing $P_{\gamma\delta}$.

To summarise, the KLV polynomials (and therefore the multiplicities  $M(\gamma,\delta)$) are completely determined by  the 
$\caH(W(\xi_a))$-module structure of  $\frM(\frB)$, and this structure is in turn completely determined by 
 the  data associated to all  $\gamma\in  \frB$
in Section \ref{dataonpar}.

The following corollary was stated and used in \cite{MAT04} and  \cite{GI}.
\begin{cor} \label{coronblocks} Suppose we have two reductive groups $G$ and $G'$ in the class of groups we consider, 
two blocks $\frB$ and $\frB'$ of Langlands-Vogan parameters with non singular infinitesimal characters, 
respectively for $G$ and $G'$, and a bijection 
\[\iota: \,  \frB\longrightarrow \frB' \]
 which respects the data associated to any $\gamma\in \frB$ (resp. $\gamma'\in \frB'$) in Section \ref{dataonpar}.
 Then 
 \[  M^G(\gamma,\delta)=M^{G'}(\iota(\gamma),\iota(\delta)), \quad (\gamma,\delta\in \frB). \]
\end{cor}

\section{Data in $G$ versus data in $L^\flat$}\label{LvsG}
Let us fix now a parabolic subgroup $P^\flat=M^\flat A^\flat N^\flat$ of $G$ 
with  $ \theta$-stable Levi factor $L^\flat=M^\flat A^\flat$.
We also fix a fundamental $\theta$-stable Catan subgroup $H^\flat$ of $L$. Such a Cartan subgroup has a decomposition
$H^\flat=T^\flat  A_{M^\flat}  A^\flat$. Of course there are similar decompositions for the Cartan subalgebras.

All the notation and results in Section \ref{Not} apply to the group $L^\flat$ instead of $G$. When needed, we will add a superscript $G$ or $L^\flat$
to distinguish between  objects defined with respect  to $G$ or  $L^\flat$.

Cartan subalgebras of $\frl_0^\flat$ are Cartan subalgebras of $\frg_0$ and Cartan subgroups of $L^\flat$ are Cartan subgroups of $G$.
This is particular the case for $ \frh_0^\flat$ and $H^\flat=T^\flat  A_{M^\flat}^\flat A^\flat$.

In general,  a Cartan subgroup $H$ in $L^\flat$  decomposes as $H=H_{M^\flat} A^\flat$ with $H_{M^\flat}$ a Cartan subgroup of 
$M^\flat$. If $H$ is $\theta$-stable, one can further decompose $H_{M^\flat} $ as $H_{M^\flat} =T A_{M^\flat}$ and $H$ as 
$H=T A_{M^\flat}A^\flat$. For such a Cartan subgroup, we have as in Section \ref{Psub}
\begin{equation} \label{comproots} R(\frg,\frh)=R(\frl^\flat,\frh)\coprod R(\frn^\flat,\frh) \coprod (-R(\frn^\flat,\frh))\end{equation}
The roots $\alpha\in R(\frn^\flat,\frh)$ are either real, or complex with $\sigma(\alpha)=-\theta(\alpha)$ also in   $R(\frn^\flat,\frh)$. 
Therefore
\[  R_{i\bbR}(\frg,\frh)=R_{i\bbR}(\frl^\flat,\frh), \quad R_{i\bbR,c}(\frg,\frh)=R_{i\bbR,c}(\frl^\flat,\frh) \]
We can therefore simply write  $R_{i\bbR}$ and $R_{i\bbR,c}$ for these systems of imaginary roots.

For the definition of the Hirai order used in the next proposition, see \cite{Hir}.
\begin{prop}
The following conjugacy classes of Cartan subgroups are in natural one-to-one correspondences.
\begin{itemize}
\item[a)] $G$-conjugacy classes of Cartan subgroups of $G$  containing a $G$-conjugate of  $A^\flat$.
\item[b)]  $K$-conjugacy classes of $\theta$-stable Cartan subgroups of $G$  containing  a $G$-conjugate of $A^\flat$.
\end{itemize}

The following conjugacy classes of Cartan subgroups are in natural one-to-one correspondences.
\begin{itemize}
\item[c)] $L^\flat$-conjugacy classes of Cartan subgroups of $L^\flat$.
\item[d)] $M_K^\flat $-conjugacy classes of  $\theta$-stable Cartan subgroups of $L^\flat$.
\item[e)] $M^\flat$-conjugacy classes of Cartan subgroups of $M^\flat$.
\item[f)] $M^\flat _K$-conjugacy classes of  $\theta$-stable Cartan subgroups of $M^\flat$.
\end{itemize}
Furthermore, the natural map from the set of conjugacy classes in $c)$ to the set of conjugacy classes in $a)$ is surjective.
The $G$-conjugacy classes in $a)$
are exactly the ones which are greater than $H^\flat$ in the Hirai order for $G$. In particular, it contains the maximally split 
$G$-conjugacy class of Cartan subgroups of $G$.
\end{prop}
\proof   The equivalence of $a)$ and $b)$  is in \cite{VGreen}, Lemma 0.1.6. Of course, it gives also the equivalence between $e)$ and $f)$ and $c)$ and $d)$.  
The equivalence of $c)$ and $e)$ is clear since $A^\flat$ is central in $L^\flat$. A Cartan subgroup containing $A^\flat$ is contained in
 $G^{A^\flat}=L^\flat$ since Cartan subgroups are abelian for linear groups, proving that  the natural map from
  the set of conjugacy classes in $c)$ to the set of conjugacy classes in $a)$ is surjective. 
  This map respects the Hirai  order (for $L^\flat$ and $G$ respectively),
   and from  this we get that the set of conjugacy classes in $a)$ are greater than the one of $H^\flat$ in the Hirai order for $G$.
   Conversely, a Cartan subgroup of $G$ with $G$-conjugacy class greater  than the one of $H^\flat$ in the Hirai order 
   for $G$ has a $G$-conjugate containing  $A^\flat$.
  If two Cartan subgroups of $L^\flat $ are $L^\flat$-conjugate, they are $G$-conjugate. In general, two  $G$-conjugate Cartan subgroups of $L^\flat$ 
 are not $L^\flat$-conjugate, unless they are maximally compact or split in $L^\flat$.
\qed

\medskip

We now fix an infinitesimal character $\chi=\chi_\xi$, both for $G$ and $L^\flat$, by choosing $\xi\in (\frh^\flat)^*$.
We decompose $\xi$ as $\xi=\xi_{M^\flat}+\nu$, according to the decomposition $(\frh^\flat)^*=(\frh_{M^\flat}^\flat)^*\oplus (\fra^\flat)^*$.

 \begin{hyp}\label{hyp1}
Consider the following conditions on $\xi=\xi_{M^\flat}+\nu \in (\frh^\flat)^*$,
\begin{itemize}
\item[ A.]  $\xi$ is non-singular for $L^\flat$, or equivalently,  $\xi_{M^\flat}$  is non-singular for $M^\flat$ {\sl i.e.}
  \[\bil{\xi_{M^\flat}}{\check \alpha}\neq 0  \text{ for all } \alpha\in R(\frl^\flat,\frh^\flat).\]
 \item[B.] \label{nugeneric0} For all $\alpha\in R(\frn^\flat,\frh^\flat)$, 
 \[
\bil{\check \alpha}{\xi}=\bil{\check \alpha}{\xi_{M^\flat}+\nu}\notin \bbZ.
\]

\end{itemize}
\end{hyp}

\begin{rmk}\label{rmkp}
If $\xi$ satisfies Hypothesis \ref{hyp1}, B., then 
\begin{equation}\label{integralroot}
R^G(\xi)= R^G(\xi_{M^\flat}+\nu)=R^{L^\flat}(\xi_{M^\flat}+\nu)=R^{M^\flat}(\xi).
 \end{equation}
Furthermore if  $\xi$ also satisfies Hypothesis \ref{hyp1}, A., then $\xi$ is non singular  also for   $G$. 
 \end{rmk}
\medskip

Our main result is 
\begin{thm} \label{main1} Suppose  $\xi=\xi_{M^\flat}+\nu \in (\frh^\flat)^*$ satisfies  Hypotheses \ref{hyp1}, A. and B.
 Let $\pi$ be an irreducible representation of $L^\flat$ of the form $\pi=\pi_{M^\flat}\boxtimes \chi_\nu$
with infinitesimal character $\chi_\xi$ .
 Then  the induced representation $i_{P^\flat}^G(\pi)= i_{P^\flat}^G(\pi_{M^\flat}\boxtimes \chi_\nu) $ is irreducible.
\end{thm}

\begin{rmk} We see that under hypothesis A., Theorem \ref{thmintro} is a corollary of the result above, since  condition B. is generic in $\nu$,
{\sl i.e} it holds for $\nu\in (\fra^\flat)^*$ outside a locally finite, countable number of affine hyperplanes. 
\end{rmk}

We will prove this theorem following the ideas given in the introduction (see corollary \ref{coronblocks}). 
We start by comparing the parameters for irreducible representations  with infinitesimal character $\chi=\chi_\xi$, for $G$ and $L^\flat$.
To be coherent with our preceding notation, we also fix a dominant  $\xi_a$ in the dual of the abstract Cartan subalgebra
$\frh_a$ such that $\chi_{\xi_a}=\chi_\xi=\chi$.

Consider a parameter $\gamma^{L^\flat}=(H=TA, \Gamma, \bar \gamma)\in \caP^{L^b}_\chi$
 as in Section  \ref{Lpar} (but for $L^\flat$). Since imaginary roots are the same for $\frl^\flat$ and $\frg$, it is clear 
that it is also a parameter in $ \caP^{G}_\chi$ (see $(2)$ and $(3)$ below) and conversely  if $H$ is greater than $H^\flat$ in the Hirai order.
 So the identity map 
\[   \caP^{L^\flat}_\chi  \longrightarrow \caP^{G}_\chi, \quad \gamma^L\mapsto \gamma^G, \] 
induces a map 
\begin{equation}\label{LGcor}
{ \caP^{L^\flat}_\chi}_ {/\sim M_K^\flat}  \longrightarrow {\caP^{G}_\chi}_{/\sim K}, \quad \gamma^L\mapsto \gamma^G, \end{equation}
with image the set of parameters $\eta=(H=TA, \Gamma, \bar \gamma)\in \caP^{G}_\chi$ with $H$  greater than $H^\flat$ in the Hirai order.

\begin{prop}\label{1to1}  Under Hypotheses \ref{hyp1}, A. and B.     the map (\ref{LGcor}) is injective. \end{prop}
\proof We have to show that if  two parameters $\gamma_1=(H_1=T_1A_1, \Gamma_1, \bar \gamma_1)$ and
 $\gamma_2=(H_2=T_2A_2, \Gamma_2, \bar \gamma_2)\in \caP^{L^b}_\chi$ are $G$-conjugate, then they are $L^\flat$-conjugate.
Since $\bar \gamma_1$ and $\bar \gamma_2$ define the same infinitesimal character as $\xi$, there are elements 
$l_1$ and $l_2$ in the complex  group $L^\flat(\bbC)$ such that $l_1\cdot \bar \gamma_1=\xi=l_2\cdot \bar \gamma_2$ and
$l_1\cdot \frh_1=\frh^\flat=l_2\cdot \frh_2$.  
Since $\gamma_1$ and $\gamma_2$ are $G$-conjugate, there is an element $g\in G$ such that $g\cdot \frh_1=\frh_2$ and
$g\cdot \bar \gamma_1=\bar \gamma_2$.
Therefore, setting $n=l_2 g l_1^{-1}\in G(\bbC)$, we get  $n\cdot \xi=\xi$ with $n\in \mathrm{Norm}_{G(\bbC)}(\frh^\flat)$. 
Since $\xi$ is non-singular, we must have  $n\in \mathrm{Centr}_{G(\bbC)}(\frh^\flat)=H^\flat\subset L^\flat(\bbC)$ and so 
$g\in L^\flat(\bbC)$. Thus   $g\in G\cap L^\flat(\bbC)=L^\flat$. \qed

\bigskip
We now check that the data   $(1)$ to $(8)$ in Section \ref{dataonpar}
associated to parameters are preserved by the correspondence $\gamma^{L^\flat}\mapsto \gamma^G$
in (\ref{LGcor}). So, let us fix    $\gamma^{L^\flat}= (H=TA, \Gamma, \bar \gamma) \in \caP^{L^\flat}_\chi$. We add superscript
$G$ or $L^\flat$ to the various objects defined in Section \ref{Lpar} to distinguished the ones defined with respect to $G$ from the ones 
defined with respect to $L^\flat$.

\begin{itemize} 
\item[(1)]  Because of hypothesis  \ref{nugeneric0}.B  on $\nu\in (\frh^\flat)^*$, we have $R^G(\xi_{M^\flat}+\nu)=R^{L^\flat}(\xi_{M^\flat}+\nu)$
and therefore $R^{L^\flat}(\bar \gamma)=R^G(\bar \gamma)$. Thus the integral root systems for $\bar \gamma$ are the same for 
$L^\flat$ and  $G$.

\item[(2) and (3)]  We have seen  that roots in $R(\frn^\flat,\frh)$ are either real, or complex with $\sigma(\alpha)=-\theta(\alpha)$. 
Therefore, the imaginary roots 
for $\frh$ are the same in $\frl^\flat$ and $\frg$, and such an imaginary root is compact in $G$ if and only if it is compact in $L^\flat$, and so we have 
\[  R^G_{i\bbR}=  R^{L^\flat }_{i\bbR}, \quad   R^G_{i\bbR,c}=  R^{L^\flat }_{i\bbR,c}\quad R^G_{i\bbR,nc}=  R^{L^\flat }_{i\bbR,nc}\]
\[  R^{G,+}_{i\bbR}=  R^{L^\flat,+ }_{i\bbR}, \quad   R^{G,+}_{i\bbR,c}=  R^{L^\flat,+ }_{i\bbR,c}\quad R^{G,+}_{i\bbR,nc}=  R^{L^\flat ,+}_{i\bbR,nc}\]
  \[\rho(R^{G,+}_{i\bbR})= \rho(R^{L^\flat,+}_{i\bbR}), \quad   \rho(R^{G,+}_{i\bbR,c})= \rho(R^{L^\flat ,+}_{i\bbR,c})\]

\item[(4)] The fact that a non compact imaginary root $\tilde \alpha$ is type I or type II in $G$ depends only 
on the map $\Phi_\alpha :\SL(2,\bbR)\longrightarrow G$ in  (\ref{Phia}) as it is clear from the equivalent conditions defining  type I or type II.
Since in our context  we can choose the same map $\Phi_\alpha :\SL(2,\bbR)\longrightarrow L^\flat \subset G$ for $G$ and $L^\flat$,  we see 
that a non compact imaginary root is type I in $L^\flat$ if and only if it is type I in $G$.

\item[(5)] Since $R^G(\bar \gamma)=R^{L^\flat}(\bar \gamma)$, we have also $R_{\bbR}^G(\bar \gamma)=R_\bbR^{L^\flat}(\bar \gamma)$
and likewise $R_{\bbR}^{G,+}(\bar \gamma)=R_\bbR^{L^\flat,+}(\bar \gamma)$.

\item[(6)] Let $\alpha \in  R^G(\bar \gamma)=R^{L^\flat}(\bar \gamma)$ be an integral real root.
We want to compare the parity condition for  $L^\flat$ and $G$.
Since $m_\alpha$ is defined via the same map $\Phi_\alpha$ for $L^\flat$ and $G$, we have only to check that 
$\epsilon_\alpha^{L^\flat}=\epsilon_\alpha^{G}$. We do that in the lemma below.
Therefore the parity condition is  the same in $L^\flat$ and $G$: 
\[ R^{L^\flat}_{\bbR,0}(\bar \gamma)= R^{L^\flat}_{\bbR,0}(\bar \gamma) \text { and  }R^{L^\flat}_{\bbR,1}(\bar \gamma)= R^{L^\flat}_{\bbR,1}(\bar \gamma) . \]

\item[(7)] As for non compact imaginary roots, real roots  are of the same type (I or II) with respect to $L^\flat$ and $G$. 

\item[(8)] Since $R^G(\bar \gamma)=R^{L^\flat}(\bar \gamma)$, we have also $R_{\bbC}^G(\bar \gamma)=R_\bbC^{L^\flat}(\bar \gamma)$
and likewise $R_{\bbC}^{G,+}(\bar \gamma)=R_\bbC^{L^\flat,+}(\bar \gamma)$, 
 $R_{\bbC,0}^{G,+}(\bar \gamma)=R_{\bbC,0}^{L^\flat,+}(\bar \gamma)$ and 
$R_{\bbC,1}^{G,+}(\bar \gamma)=R_{\bbC,1}^{L^\flat,+}(\bar \gamma)$.

\end{itemize}

\begin{lemma} \label{epsLG}
For any integral  real root $\alpha \in R_\bbR(\bar \gamma)$, $\epsilon_\alpha^{L^\flat}=\epsilon_\alpha^{G}$.
\end{lemma}
\proof The sign $\epsilon_\alpha^{G}$ is defined as $(-1)^{d+1}$ where $d$ is an integer given  by one of the definitions in \cite{VGreen},
Lemma 8.3.9. For us, the most convenient is the first one, {\sl i.e.} we take $d=d_1$ in {\sl ibid}.
We fix $\Phi_\alpha : \SL(2,\bbR)\rightarrow L^\flat$ as in Section \ref{SCT}. Thus we get $m_\alpha$ and $H^\alpha=T^\alpha A^\alpha$.
Consider a cuspidal parabolic subgroup $P_\alpha=M_\alpha A_\alpha N_\alpha$ attached to $H^\alpha$, {\sl i.e.}
$L^\alpha=M_\alpha A_\alpha=G^{A_\alpha}$. Up to conjugacy in $L^\flat$, we may assume that $A^\flat \subset A^\alpha$
since $H^\alpha$ is a Cartan subgroup of $L^\flat$, and therefore greater than $H^\flat$ in the Hirai order.
Thus $M_\alpha\subset M^\flat$.
Therefore the integer $d_1$ defined in \cite{VGreen},
Lemma 8.3.9, which is $d_1=\frac{1}{2} \dim\left( (-1)- \text{ eigenspace of  } m_\alpha \text{ in } \frmm_\alpha \cap \frk \right)$
equals 
\[\frac{1}{2} \dim\left( (-1)- \text{ eigenspace of  } m_\alpha \text{ in } \frmm_\alpha \cap \frk_{\frl^\flat} \right),\]
since $ \frmm_\alpha \cap \frk= \frmm_\alpha\cap \frmm^\flat\cap   \frk= \frmm_\alpha \cap \frk_{\frl^\flat}$.
Therefore $\epsilon_\alpha^{L^\flat}=(-1)^{d_1+1}=\epsilon_\alpha^{G}$.
\qed

\medskip 

Let us now consider the decomposition of the parameter sets ${\caP^{L^\flat}_\chi}_ {/\sim M_K^\flat} $ and $ {\caP^{G}_\chi}_{/\sim K}$ into blocks
as in (\ref{blocks}). From the characterization of blocks in terms of Cayley transforms and cross-action, we see 
that the injective correspondence $\gamma^{L^\flat}\mapsto \gamma^G$ respects blocks.

\begin{lemma} \label{surj}
Let us consider a block $\frB^G$ in  $ {\caP^{G}_\chi}_{/\sim K}$. Then all elements in $\frB^G$ are in the image of 
(\ref{LGcor}), or none of them are. Therefore, (\ref{LGcor}) induces a bijection between corresponding blocks.
\end{lemma}
\proof This is clear from the characterisation  of blocks given in \S \ref{blocks}. Indeed, suppose that  in $\frB^G$, there is a parameter $\eta=(H=TA,...)$
which is in the image of  (\ref{LGcor}) and one  $\eta'=(H'=T'A',...)$ which is not. Then $H$ is greater or equal to $H^\flat$ in the Hirai order, and 
$H'$ is not. Furthermore, there would be a sequence
of parameters $\eta_0=\eta$, $\eta_1$, ..., $\eta_r=\eta'$ in $\frB^G$ such that $\eta_{i+1}$ is obtained from $\eta_{i}$ either by the cross-action 
with respect to a Cayley transform associated to a  real or non compact imaginary simple integral root or by the cross-action of a complex 
simple integral root. Since cross-action doesn't change the conjugacy class of the associated Cartan subgroup,  there is an index $i$ such that 
 $\eta_{i}$ is in the image of (\ref{LGcor}),  $\eta_{i+1}$ is not, and furthermore $\eta_i$ and $\eta_{i+1}$ are related by 
a  Cayley transform associated to a  real or non compact imaginary simple integral root. Our problem is reduced to the case $\eta=\eta_i$ and 
$\eta'=\eta_{i+1}$.  But then $\eta'$ would also be in the image of  (\ref{LGcor}).
\qed
 \medskip

Given a block $\frB^{L^\flat}$ in 
${\caP^{L^\flat}_\chi }_{/\sim M_K^\flat} $,   we see that   the integral length (Def. \ref{li}) is the same for $\frB^{L^\flat}$ and the corresponding block 
$\frB^{G}$, if we choose the constants $c_0^G$ and $c_0^{L^\flat}$ to be equal.

\section{The case of singular infinitesimal character}\label{Sing}

We turn now to the case of possibly singular infinitesimal character $\chi=\chi_\xi$ with 
$\xi\in{\frh^\flat}^*$. The relevant discussion may be found in \cite{ABV}, Chapter 11 and \cite{V1024}.
First, the parameters have to be enriched by an extra piece of  data, so a parameter is now a multiplet 
 \[\gamma=(H=TA,\Gamma,\bar \gamma,R_{i\bbR}^+),\]
where $(H=TA,\Gamma,\bar \gamma)$ is as  before and  $R_{i\bbR}^+$ is a system of positive imaginary roots of $\frh$ in $\frg$.
The conditions imposed on $\gamma$ are the following $a)$, $b)$, $c)$,  $d)$ and $e)$ :

\begin{itemize}
\item[a)] $\bil{\alpha}{\bar \gamma}\geq  0$,  $(\forall \alpha\in R_{i\bbR}^+)$.
\end{itemize}
With $R^+_{i\bbR,c}$, $\rho(R^+_{i\bbR})$ and $ \rho(R^+_{i\bbR,c})$ as above, conditions $b)$ and $c)$ are the same as in Definition \ref{defpar}.
\begin{itemize}
\item[d)] Suppose $\alpha$ is a simple root in $R_{i\bbR}^+$ such that $\bil{\alpha}{\bar \gamma}=0$. Then $\alpha$ is non compact.
\item[e)]  Suppose $\alpha$ is a real  root in $R(\frg,\frh)$ such that $\bil{\alpha}{\bar \gamma}=0$. Then $\alpha$ does not satisfy the parity condition, 
{\sl i.e.} $\Gamma(m_\alpha)=-\epsilon_\alpha^G$.
\end{itemize}

 \begin{hyp}\label{hyp2}
Consider the following conditions on $\xi=\xi_{M^\flat}+\nu \in (\frh^\flat)^*$,
\begin{itemize}
 \item[B.] \label{nugeneric} For all $\alpha\in R(\frn^\flat,\frh^\flat)$, 
 \[
\bil{\check \alpha}{\xi}=\bil{\check \alpha}{\xi_{M^\flat}+\nu}\notin \bbZ.
\]
\item[C1.] 
For any $w\in W(\frg,\frh^\flat)$ such that $\xi_{M^\flat}-w\cdot \xi_{M^\flat}$ is non-zero, 
$\nu$ is not in the strict affine subspace in $(\fra^\flat)^*$ of solutions of  $w\cdot \nu-\nu=\xi_{M^\flat}-w\cdot \xi_{M^\flat}$.
  
 \item[C2.]   For all  $\alpha\in R(\frn^\flat,\frh^\flat)$, $\bil{\check \alpha}{\nu}\neq 0$.

\item[D.]  For  any  $ w\in W(\frg,\frh^\flat)$,  $w\cdot \xi=\xi$  implies 
$w\in W(\frl^\flat,\frh^\flat)$.
\end{itemize}
\end{hyp}

\begin{rmk} Condition B. is the same as in Hypotheses \ref{hyp1}, and we replace condition A. there, which is the 
assumption of non singular infinitesimal character, by either condition C. (meaning C1. and C2.) or condition D.
\end{rmk}

We start with the analog of Prop \ref{1to1}.
\begin{prop}\label{1to1bis}  Under Hypotheses \ref{hyp2}, C1. and C2., or Hypothesis \ref{hyp2},  D.   the map (\ref{LGcor}) is well-defined and injective.
\end{prop}
\proof We first  have to check that the extra conditions d) and e) in the definition of the parameters are preserved, but this is straightforward.
(See Lemma \ref{epsLG} for condition e)). 
 Starting the proof for injectivity as in the proof of Proposition \ref{1to1}, 
with $l_1\cdot \frh_1=\frh^\flat$, $l_1\cdot \bar \gamma_1=\xi$, $l_2\cdot \frh_2=\frh^\flat$, $l_2\cdot \bar \gamma_2=\xi$,
and $g\cdot \frh_1=\frh_2$, $g\cdot \bar \gamma_1=\bar\gamma_2$, 
we get 
$n\cdot \xi=\xi$ with $n\in \mathrm{Norm}_{G(\bbC)}(\frh^\flat)$ 
  and we conclude  under Hypothesis  \ref{hyp2},  D. as in the proof of Prop \ref{1to1}.
If we assume instead  Hypothesis  \ref{hyp2},  C. we rewrite  $n\cdot \xi=\xi$ as $\xi_{M^\flat}-n\cdot \xi_{M^\flat}=n\cdot \nu-\nu$.
So if the linear map 
\[ \phi_n: \, (\fra^\flat)^*\longrightarrow  (\frh^\flat)^*, \quad \nu\mapsto n\cdot \nu-\nu \]
is non zero, its kernel is a strict subspace of $(\fra^\flat)^*$ and the set of solutions in $(\fra^\flat)^*$ of $n\cdot \nu-\nu= \xi_{M^\flat}-n\cdot \xi_{M^\flat}$
 is strict affine subspace. Since $ \xi_{M^\flat}-n\cdot \xi_{M^\flat} $ takes only a finite number of non-zero values for 
 $n\in \mathrm{Norm}_{\bbG(\bbC)}(\frh^\flat)$,
we see that for $\nu$ outside a finite number of strict affine subspaces in $(\fra^\flat)^*$, $\xi_{M^\flat}-n\cdot \xi_{M^\flat}=n\cdot \nu-\nu$
implies $n\cdot \nu=\nu$.  Since Hypothesis   \ref{hyp2},  C2.  implies that  $\frg^\nu=\frl$ and since $L^\flat(\bbC)$ is connected 
 $G(\bbC)^\nu=L^\flat(\bbC)$. We can then conclude as in the proof of Prop \ref{1to1} that $g\in L^\flat$.
\qed

\medskip

We now use the results of \cite{ABV}, Chapter 16, using what is called there a translation datum to reduce the problem  to the case of non singular 
infinitesimal character. The translation datum consists of our singular infinitesimal character $\xi$, a weight $\mu$ for
 $H^\flat\cap M^\flat=T^\flat  A_{M^\flat}$  satisfying $\xi'=\xi+\mu$ such that 
 \begin{itemize}
\item[a)] $\xi'$ is non-singular for $\frl^\flat$
\item[b)]  If $\bil{\check \alpha}{\xi}$  is a positive integer for $\alpha\in R(\frg,\frh^\flat)$, then $\bil{\check \alpha}{\xi'}$ is a positive integer.
\end{itemize}
Set $\chi=\chi_\xi$, $\chi'=\chi_{\xi'}$.
Then by {\sl ibid}, (16.5)(a) and (16.5)(d), there is a injection $\iota^{L^\flat}_{\xi,\xi'}: \caP^{L^\flat}_\chi \rightarrow  \caP^{L^\flat}_{\chi'}$ respecting $K_M$-conjugacy classes
such that for $\gamma^{L^\flat}, \, \delta^{L^\flat}\in \caP^{L^\flat}_\chi $ : 
\[ M^{L^\flat}(\gamma^{L^\flat}, \delta^{L^\flat})=M^{L^\flat}(\iota^{l^\flat}_{\xi,\xi'}(\gamma^{L^\flat}),\iota^{l^\flat}_{\xi,\xi'}(\delta^{L^\flat})).\]

 We can use the same translation datum, but this time for $G$, and we get 
a injection $\iota^G_{\xi,\xi'}: \caP^{G}_\xi \rightarrow  \caP^{G}_{\xi'}$ respecting $K$-conjugacy classes
such that for $\gamma^{G}, \, \delta^{G}\in \caP^{G}_\xi $ : 
\[ M^{G}(\gamma^{G}, \delta^{G})=M^{G}(\iota^G_{\xi,\xi'}(\gamma^{G}),\iota_{\xi,\xi'}(\delta^{G})).\]

By the result for non singular infinitesimal character proved in the previous section, we have that 
\[M^{L^\flat}(\iota^{L^\flat}_{\xi,\xi'}(\gamma^{L^\flat}),\iota^{l^\flat}_{\xi,\xi'}(\delta^{L^\flat}))
=M^{G}(\iota^G_{\xi,\xi'}(\gamma^{G}),\iota_{\xi,\xi'}(\delta^{G})).\]
Since
\[ \xymatrix{  {\caP^{L^\flat}_\chi }_{/\sim K_{M^\flat}}   \ar[d]^{\caI_\xi} \ar[rr]^{\iota^{L^\flat}_{\xi,\xi'}}   &  
& {\caP^{L^\flat}_{\chi'}}_{/\sim K_{M^\flat}} \ar[d]^{\caI_{\xi'}}\\
 {\caP^{G}_\chi }_{/\sim K}\ar[rr]^{\iota^{G}_{\xi,\xi'}} &  &{ \caP^{G}_{\chi'}}_{/\sim K}  }  \]
is a commutative diagram, where the vertical maps are the injective  maps previously defined in (\ref{LGcor}) and denoted here $\caI_\xi$ and $\caI_{\xi'}$, 
we get :
\begin{equation} M^{L^\flat}(\gamma^{L^\flat}, \delta^{L^\flat})=M^{G}(\gamma^{G}, \delta^{G}) .\end{equation}
We also need to prove that  $M^{G}(\eta, \delta^{G}) =0$ if $\eta$ is not in the image of $\caI_\xi$.
Since we don't have the results on blocks we need readily available when the infinitesimal character is singular, we cannot apply 
Lemma \ref{surj} directly.
Assume $M^{G}(\eta, \delta^{G}) \neq 0$ and write $\eta'={\iota^{G}_{\xi,\xi'}}(\eta)$, ${\delta'}^G={\iota^{G}_{\xi,\xi'}}(\delta^G)$. 
Therefore $ M^{G}(\eta', {\delta'}^{G}))\neq 0$, and $\eta'$, ${\delta'}^{G}$ are in the same block.
 The map $\caI_{\xi'}$ is surjective on the  block which contains both $\eta'$ and
${\delta'}^{G}$ by Lemma
\ref{surj}, thus there exists $\omega' , \delta' \in {\caP^{L^\flat}_{\chi'}}_{/\sim K_{M^\flat}}$ with, $\caI_{\xi'}(\omega')=\eta'$,  $\caI_{\xi'}(\delta')={\delta'}^G$,
 and  $M^{L^\flat}(\omega', \delta')\neq 0$ by the results of the previous section.
In \cite{ABV}, Chapter 16, the maps $\iota_{\xi,\xi'}$ (called $\phi_\caT$ there) are defined
as the inverse of  partially defined bijective maps $\psi_\caT$, the domain of this map being given by the condition that
$\psi_\caT(\std(\sigma))\neq 0$ where $\psi_\caT$ is here the Zuckerman translation functor from $\xi'$ to $\xi$. 
From \cite{ABV}, Propositions 11.16 and 11.18, we see that this condition can be checked on the data associated to 
the parameter $\sigma$ in Section \ref{dataonpar}. Therefore, the domain of $\psi_\caT^{L^\flat}$ is the inverse image
of the domain of  $\psi_\caT^{G}$ by $\caI_{\xi'}$ since the map $\caI_{\xi'}$ preserves these data. We deduce that 
$\omega'$, and ${\delta'}^G$ are in the image of  $\iota^{L^\flat}_{\xi,\xi'}$, let's say from $\omega$ and $\delta$ respectively.
Now, by the commutativity of the diagram and the injectivity of the maps, we must have $\delta=\delta^{L^\flat}$ and $\caI_{\xi}(\omega)=\eta$.
Thus $\eta$ is in the image of  $\caI_{\xi}$.
\qed

\medskip

From this we deduce as before the following
\begin{thm} \label{main2} Suppose  $\xi=\xi_{M^\flat}+\nu \in (\frh^\flat)^*$ satisfies  either Hypotheses \ref{hyp2},   B. and  C., or 
Hypotheses \ref{hyp2},   B. and D.
 Let $\pi$ be an irreducible representation of $L^\flat$ of the form $\pi=\pi_{M^\flat}\boxtimes \chi_\nu$
with infinitesimal character $\chi_\xi$ .
 Then  the induced representation $i_{P^\flat}^G(\pi)= i_{P^\flat}^G(\pi_{M^\flat}\boxtimes \chi_\nu) $ is irreducible.
\end{thm}

\begin{rmk} We see that under hypothesis B and C., Theorem \ref{thmintro} is a corollary of the result above, since  conditions C1. and C2. are
 generic in $\nu$,
{\sl i.e} it holds for $\nu\in (\fra^\flat)^*$ outside a locally finite, countable number of strict affine subspaces. 
\end{rmk}

\section{An application} We explain how the results above  lead  to simplifications  in the proof of  \cite{MR6}, Theorems  5.3 and 5.4.
For background on Arthur packets, we refer to \cite{MR3}, specially in the context of classical real  groups.

Suppose that   $\mathbb G$ is a classical group  (symplectic or  special orthogonal, the case  of unitary groups
is similar but requires some adaptation in the formulation of some definitions and statements below) over $\bbR$, of rank N.
Let us denote by $\Std_{\mathrm G}$ the standard representation of the $L$-group of  $\mathbb G$
in  $\GL_N(\bbC)$ (see \cite{MR3}, \S 3.1), for instance if  $\mathbb G=\Sp_{2n}(\bbR)$, ${}^LG=\SO_{2n+1}(\bbC)\times W_\bbR$ and
$ \Std_{\mathrm G}$ is given by the inclusion of  $\SO_{2n+1}(\bbC)$ in $\GL_{2n+1}(\bbC)$.

 Let $\psi_G: \, W_\bbR \times \SL_2(\bbC)\rightarrow {}^LG$ be an Arthur parameter for  $G$, and set 
$\psi=\Std_G\circ \psi_G$,  that we see $\psi$ as a completely reducible representation  of  $W_\bbR\times \SL_2(\bbC)$.
In  \cite{MR3}, \S 4.1, we give an explicit decomposition of $\psi$ into a direct sum of 
irreducible representations and we define  good (and bad) parity for these. 
 The parameter $\psi$ is then written as $\psi=\psi_{gp}\oplus \psi_{bp}$ where $\psi_{gp}$ (resp. $\psi_{bp}$) is the part of good 
 (resp. bad) parity of $\psi$.\footnote{In \cite{MR3} and \cite{MR6}, written in french,    $\psi_{bp}$ is the  ``bonne parité'' part and 
$\psi_{mp}$ is the  ``mauvaise parité'' . } The bad parity part $\psi_{bp}$ can be further decomposed as 
$\psi_{bp}=\rho\oplus \rho^*$ for some representation $\rho$ of $ W_\bbR\times \SL_2(\bbC)$ in $\GL_{N_\rho}(\bbC)$, and $\rho^*$ is the 
contragredient of $\rho$. By the Arthur-Langlands correspondence for $\GL_{N_\rho}$, $\rho$ is the Arthur parameter of a representation,  
denoted again by  $\rho$, of $\GL_{N_\rho}(\bbR)$. The good parity part $\psi_{gp}$ is an Arthur parameter for a group $G'$ 
of the same type as $G$, but of rank $N-N_\rho$.
Furthermore $G'\times \GL_{N_\rho}(\bbR)$ is the Levi factor of a maximal parabolic subgroup $P$ of $G$.

Representations in the Arthur packet for  $G$ with parameter $\psi$ are obtained from irreducible representations  $\pi_{G'}$ in the 
Arthur packet for  $G'$ with parameter $\psi_{gp}$, as induced representations  $i_P^G(\pi_{G'}\boxtimes \rho)$.
Theorems  5.3 and 5.4. of \cite{MR6} state that these representations are indeed irreducible. In fact, as the main results of this paper show, 
 we get irreducibility  of  $i_P^G(\pi_{G'}\boxtimes \rho)$ for any representation  $\pi_{G'}$ of $G'$ and  for any  representation $\rho$ of 
 $\GL_{N_\rho}(\bbR)$,  if their   infinitesimal characters are the ones determined by $\psi$,  under some assumption on the infinitesimal character of $\rho$.
  
 Let us explain this for $\bbG=\Sp_{2N}$, the other cases being similar. In the usual coordinates, the infinitesimal character for a parameter 
 of good parity for $\bbG=\Sp_{2(N-N_\rho)}$ consists in  $N-N_\rho$ integers (up to the Weyl group action by permutation and sign changes), 
 while the infinitesimal character corresponding to $\rho$, which comes from the bad parity part, consists in $N_\rho$ complex numbers 
 which are not integers. 
 It is then obvious that the Hypotheses \ref{hyp2} B. is  satisfied, and for D., unfortunately, our hypothesis only implies that the element 
 $w$ is in the product of the Weyl groups for $\Sp_{2(N-N_\rho)}$ and $\Sp_{2N_\rho}$, rather than in the Weyl group of 
 $\Sp_{2(N-N_\rho)} \times \GL_{2N_\rho}$. If the coordinates of the infinitesimal character of $\rho$ don't contain pairs of the form $(a,-a)$
 (which is a condition easy to check starting from $\psi$), then Hypothesis D is satisfied and we get the irreducibility of 
  $i_P^G(\pi_{G'}\boxtimes \rho)$.

In general, we can do the following (see \cite{MAT04} and \cite{GI} for similar strategy): 
we  apply Corollary \ref{coronblocks} for the relevant  blocks in the groups 
$\Sp_{2(N-N_\rho)}(\bbR)\times \Sp_{2N_\rho}(\bbR)$  and $\Sp_{2N}(\bbR)$ rather than $\Sp_{2(N-N_\rho)}(\bbR)\times \GL_{N_\rho}(\bbR)$
and $\Sp_{2N}(\bbR)$, so that this time Hypothesis D is satisfied. Then, the problem is to show that the representation parabolically  induced from 
$\rho$ to $\Sp_{2N_\rho}(\bbR)$ (using the Siegel parabolic subgroup of $\Sp_{2N_\rho}(\bbR)$) is irreducible.
This is a particular case of our original problem, but the arguments in the proof in \cite{MR6} are then technically easier. 
Other approaches to this problem may work, for instance  $\rho$ being  unitary, one may start by using Tadic's classification of the unitary dual
of general linear groups to write it in terms of Speh representations (which can be done directly from $\psi_{bp}$) and then try to use
the independence of ``polarization results'' of \cite{KV}, Chapter XI, to reduce further the problem.

\bibliographystyle{plainurl}
\bibliography{IRR_IND_D}

\begin{thebibliography}{10}

\bibitem{ABV}
Jeffrey Adams, Dan Barbasch, and David~A. Vogan, Jr.
\newblock {\em The {L}anglands classification and irreducible characters for
  real reductive groups}, volume 104 of {\em Progress in Mathematics}.
\newblock Birkh\"auser Boston, Inc., Boston, MA, 1992.
\newblock URL: \url{http://dx.doi.org/10.1007/978-1-4612-0383-4}, \href
  {http://dx.doi.org/10.1007/978-1-4612-0383-4}
  {\path{doi:10.1007/978-1-4612-0383-4}}.

\bibitem{Dat}
J.-F. Dat.
\newblock Repr\'{e}sentations lisses {$p$}-temp\'{e}r\'{e}es des groupes
  {$p$}-adiques.
\newblock {\em Amer. J. Math.}, 131(1):227--255, 2009.
\newblock URL: \url{https://doi.org/10.1353/ajm.0.0041}, \href
  {http://dx.doi.org/10.1353/ajm.0.0041} {\path{doi:10.1353/ajm.0.0041}}.

\bibitem{GI}
Wee~Teck Gan and Atsushi Ichino.
\newblock On the irreducibility of some induced representations of real
  reductive {L}ie groups.
\newblock {\em Tunis. J. Math.}, 1(1):73--107, 2019.
\newblock URL: \url{https://doi.org/10.2140/tunis.2019.1.73}, \href
  {http://dx.doi.org/10.2140/tunis.2019.1.73}
  {\path{doi:10.2140/tunis.2019.1.73}}.

\bibitem{Hir}
Takeshi Hirai.
\newblock Invariant eigendistributions of {L}aplace operators on real simple
  {L}ie groups. {II}. {G}eneral theory for semisimple {L}ie groups.
\newblock {\em Japan. J. Math. (N.S.)}, 2(1):27--89, 1976.
\newblock URL: \url{https://doi.org/10.4099/math1924.2.27}, \href
  {http://dx.doi.org/10.4099/math1924.2.27} {\path{doi:10.4099/math1924.2.27}}.

\bibitem{KV}
Anthony~W. Knapp and David~A. Vogan, Jr.
\newblock {\em Cohomological induction and unitary representations}, volume~45
  of {\em Princeton Mathematical Series}.
\newblock Princeton University Press, Princeton, NJ, 1995.

\bibitem{MAT04}
Hisayosi Matumoto.
\newblock On the representations of {${\rm Sp}(p,q)$} and {${\rm SO}^*(2n)$}
  unitarily induced from derived functor modules.
\newblock {\em Compos. Math.}, 140(4):1059--1096, 2004.
\newblock URL: \url{https://doi.org/10.1112/S0010437X03000629}, \href
  {http://dx.doi.org/10.1112/S0010437X03000629}
  {\path{doi:10.1112/S0010437X03000629}}.

\bibitem{MR6}
Colette M{\oe}glin and David Renard.
\newblock Sur les paquets d'arthur des groupes unitaires et quelques
  conséquences pour les groupes classiques.
\newblock {\em Pac. J. Math.}, 299(1):53--88, 2019.
\newblock \href {http://dx.doi.org/10.2140/pjm.2019.299.53}
  {\path{doi:10.2140/pjm.2019.299.53}}.

\bibitem{MR3}
Colette M{\oe}glin and David Renard.
\newblock Sur les paquets d'arthur des groupes classiques réels.
\newblock {\em J. Eur. Math. Soc. (JEMS)}, 22(6):1827--1892, 2020.
\newblock \href {http://dx.doi.org/10.4171/JEMS/957}
  {\path{doi:10.4171/JEMS/957}}.

\bibitem{MOY}
Matringe Nadir, Omer Offen, and Chang Yang.
\newblock On local intertwining periods.
\newblock {\em arXiv}, 2023.
\newblock URL: \url{https://arxiv.org/abs/2303.03663}, \href
  {http://dx.doi.org/https://doi.org/10.48550/arXiv.2303.03663}
  {\path{doi:https://doi.org/10.48550/arXiv.2303.03663}}.

\bibitem{Sau}
Fran\c{c}ois Sauvageot.
\newblock Principe de densit\'{e} pour les groupes r\'{e}ductifs.
\newblock {\em Compositio Math.}, 108(2):151--184, 1997.
\newblock URL: \url{https://doi.org/10.1023/A:1000216412619}, \href
  {http://dx.doi.org/10.1023/A:1000216412619}
  {\path{doi:10.1023/A:1000216412619}}.

\bibitem{VGreen}
David~A. Vogan, Jr.
\newblock {\em Representations of real reductive {L}ie groups}, volume~15 of
  {\em Progress in Mathematics}.
\newblock Birkh\"auser, Boston, Mass., 1981.

\bibitem{VIC4}
David~A. Vogan, Jr.
\newblock Irreducible characters of semisimple {L}ie groups. {IV}.
  {C}haracter-multiplicity duality.
\newblock {\em Duke Math. J.}, 49(4):943--1073, 1982.
\newblock URL: \url{http://projecteuclid.org/euclid.dmj/1077315538}.

\bibitem{VIC3}
David~A. Vogan, Jr.
\newblock Irreducible characters of semisimple {L}ie groups. {III}. {P}roof of
  {K}azhdan-{L}usztig conjecture in the integral case.
\newblock {\em Invent. Math.}, 71(2):381--417, 1983.
\newblock URL: \url{http://dx.doi.org/10.1007/BF01389104}, \href
  {http://dx.doi.org/10.1007/BF01389104} {\path{doi:10.1007/BF01389104}}.

\bibitem{V1024}
David~A. Vogan, Jr.
\newblock Understanding the unitary dual.
\newblock Lie group representations {I}, {Proc}. {Spec}. {Year}, {Univ}. {Md}.,
  {College} {Park} 1982-83, {Lect}. {Notes} {Math}. 1024, 264-286 (1983).,
  1983.

\end{thebibliography}

\end{document}